\documentclass[12pt]{article}
\usepackage{amsmath,amsfonts,amssymb,amsthm}
\newtheorem*{theor}{Theorem}
\newtheorem{theo}{Theorem}
\newtheorem{conj}{Conjecture}
\newtheorem*{conjj}{Conjecture}
\newtheorem{lem}{Lemma}
\numberwithin{equation}{section}
\newcommand\la{\lambda}

\title{Two positivity conjectures for Kerov polynomials}
\author{Michel Lassalle\\
\small Centre National de la Recherche Scientifique\\[-0.8ex]
\small Institut Gaspard-Monge, Universit\'e de Marne-la-Vall\'ee\\[-0.8ex]
\small 77454 Marne-la-Vall\'ee Cedex, France\\[-0.8ex]
\small \texttt{lassalle @ univ-mlv.fr}\\[-0.8ex]
\small \texttt{http://igm.univ-mlv.fr/{\textasciitilde}lassalle}}
\date{}

\begin{document}
\maketitle
\begin{abstract}
Kerov polynomials express the normalized characters of irreducible representations of the symmetric group, evaluated on a cycle, as polynomials in the ``free cumulants'' of the associated Young diagram. We present two positivity conjectures for their coefficients. The latter are stronger than the positivity conjecture of Kerov-Biane, recently proved by F\'eray.
\end{abstract}

\section{Kerov polynomials}

\subsection{Characters}

A partition $\la= (\la_1,...,\la_r)$ is a finite weakly decreasing
sequence of nonnegative integers, called parts. The number
$l(\la)$ of positive parts is called the length of
$\la$, and $|\la| = \sum_{i = 1}^{r} \la_i$
the weight of $\la$. For any integer $i\geq1$,
$m_i(\la) = \textrm{card} \{j: \la_j  = i\}$
is the multiplicity of the part $i$ in $\la$.

Let $n$ be a fixed positive integer and $S_n$ the group of permutations of $n$ letters. Each permutation $\sigma \in S_n$ factorizes uniquely as a product of disjoint cycles, whose respective lengths are ordered such as to form a partition $\mu=(\mu_1,\ldots,\mu_r)$ with weight $n$, the so-called cycle-type of $\sigma$. 

The irreducible representations of $S_n$ and their corresponding characters are also labelled by partitions $\la$ with weight $|\la|=n$. We write $\textrm{dim}\,\la$ for the dimension of the representation $\la$ and ${\chi}^\la_\mu$ for the value of the character $\chi^\la(\sigma)$ at any permutation $\sigma$ of cycle-type $\mu$.

Let $r\le n$ be a positive integer and  $\mu=(r,1^{n-r})$ the corresponding $r$-cycle in $S_n$. We write 
\[\hat{\chi}^\la_r=n(n-1) \cdots (n-r+1) \frac{\chi^\la_{r,1^{n-r}}}{\textrm{dim}\,\la}\] for the value at $\mu$ of the normalized character.

It was first observed by Kerov\cite{K} and Biane\cite{B2} that $\hat{\chi}^\la_r$ may be written as a polynomial in the ``free cumulants'' of the Young diagram of $\la$.

\subsection{Free cumulants}

Two increasing sequences $y=(y_1,\ldots,y_{d-1})$ and $x=(x_1,\ldots,x_{d-1},x_d)$ are said to be interlacing if $x_1<y_1<x_2<\cdots<x_{d-1}<y_{d-1}<x_d$. 
The center of the pair is $c(x,y)=\sum_i x_i-\sum_i y_i$.

To any pair of interlacing sequences with center $0$ we associate the rational function
\[G_{x,y}(z)= \frac{1}{z-x_d}\prod_{i=1}^{d-1}\frac{z-y_i}{z-x_i},\]
and the formal power series inverse to $G_{x,y}$ for composition, 
\[G_{x,y}^{(-1)}(z) = z^{-1} + \sum_{k\ge1} R_k(x,y) z^{k-1}.\] 
Note that $R_1(x,y)=c(x,y)=0$. The quantities $R_k(x,y), k\ge 2$ are called the free cumulants of the interlacing pair $(x,y)$.

Being given a partition $\la$, we consider the collection of unit boxes centered on the nodes
$\{ (j-1/2,i-1/2) : 1 \le i \, \le l(\la), 1 \le j \le {\la}_{i} \}$.
This defines a compact region in $\mathsf{R}^2$, the so-called Young diagram of $\la$. On $\mathsf{R}^2$ we define the content function by $c(u,v)=u-v$. By convention, the content of a box is the one of its center.

Then it is easily shown that the Young diagram of $\la$ defines a pair of interlacing sequences, formed by the contents $y_1,\ldots, y_{d-1}$ of its corner boxes, and the contents $x_1,\ldots,x_{d-1},x_d$ of the corner boxes of its compliment in $\mathsf{R}^2$. We have $x_1=-l(\la)$, and $x_d=\la_1$.

Conversely, every pair of interlacing sequences with integer entries and center zero uniquely determines the Young diagram of a partition $\la$. 

The free cumulants $R_k(\la), k\ge 2$ are defined accordingly. These quantities arise in the asymptotic study of representations of symmetric groups~\cite{B1}.

\subsection{Known results}

The following result was first proved in~\cite{B2} and attributed to Kerov\cite{K}.

\begin{theor}
There exist polynomials $K_r, r\ge 2$ such that for any partition $\la$ with $|\la|\ge r$, one has
\[\hat{\chi}^\la_r=K_r(R_2(\la),R_3(\la),\ldots,R_{r+1}(\la)).\]
These polynomials have integer coefficients.
\end{theor}

Let $(R_2,\ldots,R_{r+1})$ be the indeterminates of the ``Kerov polynomial'' $K_r$ and define $|\mu|$ as the ``weight'' of the monomial $R_\mu=\prod_{i\ge 2} R_i^{m_i(\mu)}$. We may decompose $K_r$ in its graded components with respect to the weight, writing 
\[K_r= \sum_{s\ge 2}  K_{r,s} \quad \mathrm{with}\quad
K_{r,s}= \sum_{|\mu|=s} x_\mu^{(r)} \prod_{i\ge 2} R_i^{m_i(\mu)}.\] 
Then it may be proved~\cite{B2} that the term of highest weight is $R_{r+1}$ and that $K_{r,s}=0$ when $s=r-2k$.

Goulden and Rattan~\cite{GR,R} have given a general formula for $K_{r,r-2k+1}$, expressing it as some coefficient in a formal power series (see also~\cite{B3}). As a consequence, one has
\[K_{r,r-1}=\frac{1}{4} \binom{r+1}{3} \sum_{|\mu|=r-1}  l(\mu)! \, \prod_{i\ge 2} \frac{((i-1)R_i)^{m_i(\mu)}}{m_i(\mu)!},\]
which had been conjectured by Biane~\cite{B2} and differently proved by \.{S}niady~\cite{S}.

The same method provides an explicit form for $K_{r,r-3}$. But as far as $K_{r,r-5}$ (and lower components) are concerned, it seems very difficult to apply. Rattan~\cite[Theorem 3.5.12]{R} found a messy expression of $K_{r,r-5}$ giving an idea about the complexity of the problem.

The following positivity property had been conjectured by Kerov~\cite{K} and Biane~\cite{B2} and was recently proved by F\'eray~\cite{F}.
\begin{theor}
The coefficients of $K_r$ are nonnegative integers.
\end{theor}

The purpose of this note is to present a stronger conjectural property.

\section{Conjectures}

An algebraic basis of the (abstract) symmetric algebra with real coefficients is formed by the classical symmetric functions, elementary $e_i$, complete $h_i$ or power-sum $p_i$. As usual for any partition $\mu$, denote $e_\mu$, $h_\mu$ or $p_\mu$ their product over the parts of $\mu$, and $m_\mu$ the monomial symmetric function, sum of all distinct monomials whose exponent is a permutation of $\mu$.

For a clearer display we write
\[\mathcal{R}_\mu=\prod_{i\ge 2} 
((i-1)R_i)^{m_i(\mu)}/m_i(\mu)!.\]

Firstly we conjecture that the Kerov components $K_{r,r-2k+1}$ may be described in a unified way, \textit{independent of} $r$.
\begin{conj}
For any $k\ge 1$ there exists an inhomogeneous symmetric function $f_k$, having maximal degree $4(k-1)$, such that
\[K_{r,r-2k+1}= \binom{r+1}{3} \sum_{|\mu|=r-2k+1} 
(l(\mu)+2k-2)!\, f_k(\mu)  \,  \mathcal{R}_\mu,\]
where $f_k(\mu)$ denotes the value of $f_k$ at the integral vector $\mu$. This symmetric function is independent of $r$.
\end{conj}

The assertion is trivial for $k=1$ since we have $f_1=1/4$. Secondly we conjecture the symmetric function $f_k$ to be positive in the following sense.
\begin{conj}
For $k\ge 2$ the inhomogeneous symmetric function $f_k$ may be written
\[f_k= \sum_{|\rho| \le 4(k-1)} c_\rho^{(k)} m_\rho,\]
where the coefficients $c_\rho^{(k)}$ are positive rational numbers.
\end{conj}

The positivity of the coefficients of $K_{r,r-2k+1}$ is an obvious consequence. We emphasize that the coefficients of $f_k$ in terms of any other classical basis \textit{may be negative}.

Conjecture 2 is firstly supported by the case $k=2$. Using the expression of $K_{r,r-3}$ given in~\cite{GR}, we have the following result, whose proof is postponed to Section 3.

\begin{theo}
For $k=2$, we have
\begin{multline*}
5760 f_2=
3 m_{4}+8 m_{31}+10 m_{22}+16 m_{21^2}+24 m_{1^4}\\
+20 m_{3}+36 m_{21}+48 m_{1^3}+35 m_{2}+40 m_{1^2}+18 m_{1}.
\end{multline*}
\end{theo}

Conjecture 2 is secondly supported by extensive computer calculations, giving the values of the positive numbers $c_\rho^{(k)}$ for $k=3, 4$. The two following conjectures have been checked for any $K_r$ with $r\le 32$.

\begin{conj}
For $k=3$, the values of $2 . 6! . 8 ! \, c_\rho^{(3)}$ are given by the table below.
\end{conj}
\small

\hspace{1 cm}\\
\begin{tabular}{|c|c|c|c|c|c|c|c|c|c|c|}
\hline
8 & 71 & 62 & $61^2$ & 53 & 521 & $51^3$ & $4^2$ & 431 & $42^2$ & $421^2$ \\ \hline
9 & 48 & 132 & 224 & 240 & 544 & 908 & 294 & 848 &1132 & 1904 \\ \hline
$41^4$ & $3^22$ & $3^21^2$ & $32^21$ & $321^3$ & $31^5$ & $2^4$ & $2^31^2$ & $2^21^4$ & $21^6$ & $1^8$ \\ \hline
3148 & 1440 & 2440 & 3280 & 5480 & 9040 & 4440 & 7440 & 12360 & 20400 & 33600\\ \hline
\end{tabular}

\hspace{0.3 cm}\\
\begin{tabular}{|c|c|c|c|c|c|c|c|c|c|c|}
\hline
7 & 61 & 52 & $51^2$ & 43 & 421 & $41^3$ & $3^21$ & $32^2$ & $321^2$ & $31^4$ \\ \hline
216 & 968 & 2296 & 3744 & 3560 & 7704 & 12368 & 9856 & 13072 & 21264 & 33968 \\ \hline
$2^31$ & $2^21^3$ & $21^5$ & $1^7$\\ 
\cline{1-4}
28560 & 46080& 73680 & 117600\\
\cline{1-4}
\end{tabular}

\hspace{0.3 cm}\\
\begin{tabular}{|c|c|c|c|c|c|c|c|c|c|c|}
\hline
6 & 51 & 42 & $41^2$ & $3^2$ & 321 & $31^3$ & $2^3$ & $2^21^2$ & $21^4$ & $1^6$ \\ \hline
2094 & 7696 & 15450 & 24016 & 19696 & 40592 & 62428 & 53796 & 83848 & 128988 & 198120 \\ \hline
\end{tabular}

\hspace{0.3 cm}\\
\begin{tabular}{|c|c|c|c|c|c|c|}
\hline
5 & 41 & 32 & $31^2$ & $2^21$ & $21^3$ & $1^5$ \\ \hline
10588 & 30972 & 51096 & 75232 & 99640 & 146200 & 214040 \\ \hline
\end{tabular}

\hspace{0.3 cm}\\
\begin{tabular}{|c|c|c|c|c|}
\hline
4 & 31 & $2^2$ & $21^2$ & $1^4$ \\ \hline
30109 & 67360 & 87382 & 120912 & 166320 \\ \hline
\end{tabular}

\hspace{0.3 cm}\\
\begin{tabular}{|c|c|c|}
\hline
3 & 21 & $1^3$ \\ \hline
48092 & 77684 & 98016 \\ \hline
\end{tabular}

\hspace{0.3 cm}\\
\begin{tabular}{|c|c|c|}
\hline
2 & $1^2$ & 1 \\ \hline
39884 & 43928 & 13200 \\ \hline
\end{tabular}

\normalsize
\newpage
\begin{conj}
For $k=4$, the values of $2.8!.12! \, c_\rho^{(4)}$ are given by the table below.
\end{conj}
\small

\hspace{1 cm}\\
\begin{tabular}{|c|c|c|c|c|c|c|}
\hline
12 & 11,1 & 10,2 & $10,1^2$ & 93 & 921 & $91^3$ \\ \hline
495 & 3960 & 16830 & 29040 & 48312 & 113520 & 194392 \\ \hline
$84$ & 831 & $82^2$ & $821^2$ & $81^4$ & 75 & 741 \\ \hline
99297 & 296472 & 403590 &  692912 & 1180248 & 150480 & 546480 \\ \hline
732 & $731^2$ & $72^21$ & $721^3$ & $71^5$ & 66 & 651 \\ \hline
945120 & 1626592 & 2219360 & 3792480 &6439200 & 172260 & 733920 \\ \hline
642 & $641^2$ & $63^2$ & 6321 & $631^3$ & $62^3$ & $62^21^2$ \\ \hline
1543740 & 2654432 & 1960992 & 4611552 & 7890528 & 6305640 & 10797440\\ \hline
$621^4$ & $61^6$ & $5^22$ & $5^21^2$ & 543 & 5421 & $541^3$ \\ \hline
18388320 & 31168800 & 1811040  & 3110800 & 2797872 & 6566560 & 11221392 \\ \hline
$53^21$ & $532^2$ & $5321^2$ & $531^4$ & $52^31$ & $52^21^3$ & $521^5$ \\ \hline
8360352 & 11420640 & 19573792 &  33343968 & 26812800 & 45753120 &  77733600\\ \hline
$51^7$ & $4^3$ & $4^231$ & $4^22^2$ & $4^221^2$ & $4^21^4$ & $43^22$\\ \hline
131644800 &  3402630 &  10157840 & 13861540 &  23751840 & 40429200 &  17629920 \\ \hline
$43^21^2$ & $432^21$ & $4321^3$ & $431^5$ & $42^4$ & $42^31^2$ & $42^21^4$ \\ \hline  
30274720 &  41416480 &  70724640 &  120150240 &  56773080 & 97050240 & 165207840 \\ \hline
$421^6$ & $41^8$ & $3^4$ & $3^321$ & $3^31^3$ & $3^22^3$ & $3^22^21^2$ \\ \hline
280274400 & 474445440 &  22397760 &  52718400 &  90162240 &  72246720 &  123618880 \\ \hline
$3^221^4$ & $3^21^6$ & $32^41$ & $32^31^3$ & $32^21^5$ & $321^7$ & $31^9$ \\ \hline
210584640 & 357315840 & 169727040 &  289477440 &  492072000 & 834301440 & 1412328960\\  \hline
$2^6$ & $2^51^2$ & $2^41^4$ & $2^31^6$ & $2^21^8$ & $21^{10}$ & $1^{12}$\\\hline
233226000 & 398160000 & 677678400 & 1150632000 &1950278400 &3302208000 & 5588352000\\\hline
\end{tabular}

\hspace{.3 cm}\\
\begin{tabular}{|c|c|c|c|c|c|c|}
\hline
11 & 10,1 & 92 & $91^2$ & 83 & 821 & $81^3$ \\ \hline
25740 & 184140 & 708444 & 1199440 & 1836252 & 4210844 & 7075728\\ \hline
$74$ & 731 & $72^2$ & $721^2$ & $71^4$ & 65 & 641 \\ \hline
3371544&9817984&13294160&22416768&37488576&4518360&15961880\\ \hline
632 & $631^2$ & $62^21$ & $621^3$ & $61^5$ & $5^21$ & 542 \\ \hline
27441744&46345728&62955728&105656256&176236800&18678880&38954344\\ \hline
$541^2$ & $53^2$ & 5321 & $531^3$ & $52^3$ & $52^21^2$ & $521^4$ \\ \hline
65688480&49454592&113453824&190456128&154321200&259558848&434150016\\ \hline
$51^6$ & $4^23$ & $4^221$ & $4^21^3$ & $43^21$ & $432^2$ & $4321^2$ \\ \hline
723211200&59989160&137408040&230425440&174697600&237252400&399329280\\ \hline
$431^4$ & $42^31$ & $42^21^3$ & $421^5$ & $41^7$ & $3^32$ & $3^31^2$
\\ \hline
667719360&544244400&912287040&1522644480&2534616000&301150080&507776640\\ \hline
$3^22^21$ & $3^221^3$ & $3^21^5$ & $32^4$ & $32^31^2$ & $32^21^4$ & $321^6$  \\ \hline
691290880&1159603200&1935373440&942671520&1583527680&2648849280&4416612480\\ \hline
$31^8$ & $2^51$ & $2^41^3$ & $2^31^5$ & $2^21^7$ & $21^9$ & $1^{11}$  \\ \hline
7350658560&2163722400&3624808320&6054048000&10089797760 &16795537920 & 27941760000\\ \hline
\end{tabular}

\hspace{.3 cm}\\
\begin{tabular}{|c|c|c|c|c|c|}
\hline
10 & 91 & 82 & $81^2$ & 73 & 721  \\ \hline
589545&3732696&12880197&21347832&29796624&66542608\\ \hline
$71^3$ & $64$ & 631 & $62^2$ & $621^2$ & $61^4$   \\ \hline
109503504&48249234&136592720&184006988&304004800&498221712\\ \hline
$5^2$ & 541 & 532 & $531^2$ & $52^21$ & $521^3$ \\ \hline
56379312&192905680&329380304&544358320&
736055232&1210234416\\ \hline
$51^5$ & $4^22$ & $4^21^2$ & $43^2$ & $4321$ & $431^3$ \\ \hline
1979174160&398071454&657018384&504522128&1126245296&1850729904\\ \hline
$42^3$ & $42^21^2$ & $421^4$ & $41^6$ & $3^31$ & $3^22^2$ \\ \hline
1523067348&2510092224&4113855024&6719636880&1427727840&1928190880
\\ \hline
$3^221^2$ & $3^21^4$ & $32^31$ & $32^21^3$ & $321^5$ & $31^7$ 
\\ \hline
3180030560&5210415840&4309828320&7080806880&11585384160&18913547520 \\ \hline
$2^5$ & $2^41^2$ & $2^31^4$ & $2^21^6$ & $21^8$ & $1^{10}$ \\ 
\hline
5844598200&9618960960&15768879840&25780980960&42087911040&68660524800\\ \hline
\end{tabular}

\hspace{.3 cm}\\
\begin{tabular}{|c|c|c|c|c|c|}
\hline
9 & 81 & 72 & $71^2$ & 63 & 621 \\ \hline
7834926&43370910&132689304&214757664&270145656&585908840\\ \hline
$61^3$ & $54$ & 531 & $52^2$ & $521^2$ & $51^4$ \\ \hline
942097728&380072484&1041283232&1395178488&2251722880&3607638624\\ \hline
$4^21$ & 432 & $431^2$ & $42^21$ & $421^3$ & $41^5$ \\ \hline
1253522292&2121348680&3421048224&4600109272&7388286912&11813196960\\ \hline
$3^3$ & $3^221$ & $3^21^3$ & $32^3$ &$32^21^2$ & $321^4$\\ \hline
2679266304&5805122752&9318556608&7798935408&12559063744&20114667264\\ \hline
$31^6$ & $2^41$ & $2^31^3$ & $2^21^5$ & $21^7$ & $1^9$ \\ \hline
32123903040&16915888080&27152536320&43428598080&69327800640&110563004160\\ \hline
\end{tabular}

\hspace{.3 cm}\\
\begin{tabular}{|c|c|c|c|c|}
\hline
8 & 71 & 62 & $61^2$ & 53 \\ \hline
66992805&319460328&854070228&1345992736&1504935432\\ \hline
521 & $51^3$ & $4^2$ & 431 & $42^2$ \\ \hline
3156966208&4945126296&1806665454&4760982424&6336879340\\ \hline
$421^2$ & $41^4$ & $3^22$ & $3^21^2$ & $32^21$ \\ \hline
9953455776&15535885752&7959879312&12492469616&16671548080\\ \hline
$321^3$ & $31^5$ & $2^4$ & $2^31^2$ & $2^21^4$ \\ \hline
26065233552&40592042160&22229994072&34840460832&54337307568\\ \hline
$21^6$ & $1^8$ \\ \cline{1-2}
84517248240& 131257445760  \\ \cline{1-2}
\end{tabular}

\hspace{.3 cm}\\
\begin{tabular}{|c|c|c|c|c|}
\hline
7 & 61 & 52 & $51^2$ & 43 \\ \hline
386137224&1557181296&3572220960&5460878192&5341858632\\ \hline
421 & $41^3$ & $3^21$ & $32^2$ & $321^2$ \\ \hline
10769122320&16360041456&13438992512&17722898864&26967001248
\\ \hline
$31^4$ & $2^31$ & $2^21^3$ & $21^5$ & $1^7$\\ \hline
40796325216&35619645600&53958337440&81409500480&122509104480
\\ \hline
\end{tabular}

\hspace{.3 cm}\\
\begin{tabular}{|c|c|c|c|c|c|}
\hline
6 & 51 & 42 & $41^2$ & $3^2$ & 321 \\ \hline
1527234687&5086528128&9789272361&14430109232&12134469600&23282303088\\ \hline
$31^3$ & $2^3$ & $2^21^2$ & $21^4$ & $1^6$ \\ \cline{1-5}
34060600640&30307366254&44384647296&64583789280&93548535360
 \\ \cline{1-5}
\end{tabular}

\hspace{.3 cm}\\
\begin{tabular}{|c|c|c|c|c|}
\hline
5 & 41 & 32 & $31^2$ & $2^21$ \\ \hline
4133557494&11019741678&17318813292&24369700608&31165644708\\ \hline
$21^3$ & $1^5$ \\ \cline{1-2}
43403668704&59946923520\\ \cline{1-2}
\end{tabular}

\hspace{.3 cm}\\
\begin{tabular}{|c|c|c|c|c|}
\hline
4 & 31 & 22 & $21^2$ & $1^4$ \\ \hline
7478442180&15298473960&19094031000&25180566840&32685206400\\ \hline
\end{tabular}

\hspace{.3 cm}\\
\begin{tabular}{|c|c|c|}
\hline
3 & 21 & $1^3$  \\ \hline
8579601096&12733485336&15147277200\\ \hline
\end{tabular}

\hspace{.3 cm}\\
\begin{tabular}{|c|c|c|}
\hline
 2 & $1^2$ & 1 \\ \hline
 5589321408 & 5773242816 & 1555424640 \\ \hline
\end{tabular}

\normalsize
\section{Proof of Theorem 1}

Following~\cite{GR,R} we consider the generating series
\[C(z)=\sum_{i\ge 0}C_iz^i= \Big(1-\sum_{i\ge 2}(i-1)R_i z^i\Big)^{-1}.\]
By classical methods we have
\begin{equation*}
C_n= \sum_{|\mu|=n} l(\mu)! \, \mathcal{R}_\mu.
\end{equation*}
It may be shown (see a proof in Section 7 below) that if $\phi$ is a polynomial in $i$, there exists a symmetric function $\hat{\phi}$ such that
\[\sum_{\begin{subarray}{c}(i,j,k) \in \mathsf{N}^3\\i+j+k=n
\end{subarray}} \phi(i)\,C_iC_jC_k = \sum_{|\mu|=n} (l(\mu)+2)! \, \hat{\phi}(\mu)\,\mathcal{R}_\mu,\]
where $\hat{\phi}(\mu)$ denotes the value of $\hat{\phi}$ at the integral vector $\mu$. For $\phi(i)=a+bi+ci^2$, we have
\[\hat{\phi}=a/2+bn/6+c(n^2+p_2)/12.\]

The following explicit form of $K_{r,r-3}$ was given in~\cite[Theorem 3.3]{GR}
\[K_{r,r-3}=\binom{r+1}{3} \sum_{\begin{subarray}{c}(i,j,k) \in \mathsf{N}^3\\i+j+k=r-3\end{subarray}}
(a(r) +b(r) i^2)C_iC_jC_k,\]
with
\[a(r)=-\frac{1}{2880}(r-1)(r-3)(r^2-4r-6),\quad\quad\quad b(r)=\frac{1}{480}(2r^2-3).\]
As a straightforward consequence, we have
\[K_{r,r-3}=\binom{r+1}{3} \sum_{|\mu|=r-3} (l(\mu)+2)! \,f_2(\mu) \, \mathcal{R}_\mu,\]
with
\[f_2(\mu)=\frac{1}{2}a(r)+\frac{1}{12}b(r)\big((r-3)^2+p_2(\mu)\big).\]
But since $|\mu|=p_1(\mu)=r-3$, this can be rewritten
\[f_2=\frac{1}{5760}\Big(2p_2p_1^2+p_1^4+12p_2p_1+8p_1^3+15p_2+20p_1^2+18p_1\Big).\]
Using for instance ACE~\cite{V} we easily obtain
\begin{multline*}
f_2=\frac{1}{5760}
\Big(3 m_{4}+8 m_{31}+10 m_{22}+16 m_{21^2}+24 m_{1^4}\\
+20 m_{3}+36 m_{21}+48 m_{1^3}+35 m_{2}+40 m_{1^2}+18 m_{1}\Big).
 \qed \end{multline*}

Observe that in this particular situation, the coefficients of $f_2$ in terms of power sums are nonnegative. This property \textit{is not true} for $K_{r,r-5}$ and lower components.

Starting from~\cite[Theorem 3.5.12]{R}, Conjecture 3 may probably be proved along the same line.

\section{$C$-expansion}

Goulden and Rattan~\cite{GR,R} have considered the expansion of Kerov polynomials in terms of the indeterminates $C_i$. They have given the following positivity conjecture, proved for $k=1,2$, which is stronger than the one of Kerov and Biane.
\begin{conjj}
For $k\ge 1$ the coefficients  of $K_{r,r-2k+1}$ in terms of the $C_i$'s are nonnegative rational numbers.
\end{conjj}

In analogy with Section 2 we  conjecture that for any $k\ge 1$ one has
\[K_{r,r-2k+1}= \binom{r+1}{3} \sum_{\begin{subarray}{c}\nu \in \mathsf{N}^{2k-1}\\|\nu|=r-2k+1\end{subarray}} F_k(\nu)
\,  \prod_{i=1}^{2k-1}C_{\nu_i},\]
where $F_k$ is an inhomogeneous symmetric function, having maximal degree $4(k-1)$ and independent of $r$.

This is clear for $k=1$ since 
\[K_{r,r-1}= \frac{1}{4}\binom{r+1}{3}C_{r-1},\]
hence $F_1=1/4$. For $k=2$ we have seen in Section 3 that 
\[F_2(\nu)=a(r)+\frac{1}{3}b(r)p_2(\nu)\]
with $\nu =(i,j,k)$. Since $|\nu|=p_1(\nu)=r-3$, we obtain
\[F_2=\frac{1}{2880}\Big(4p_2p_1^2-p_1^4+24p_2p_1-4p_1^3+30p_2+5p_1^2+18p_1\Big).\]

However we emphasize that, unlike those of $f_2$, the coefficients of $F_2$ in terms of monomial symmetric functions \textit{are not positive}. One has
\[F_2=\frac{1}{2880}
\Big(3 m_{4}+4 m_{31}+2 m_{22}-4 m_{21^2}
+20 m_{3}+12 m_{21}-24 m_{1^3}+35 m_{2}+10 m_{1^2}+18 m_{1}\Big).
\]
Therefore it seems that $C$-positivity and $R$-positivity are of a different nature.

\section{New expansion}

For a better understanding of the difference between the $C$ and $R$ expansions, it is useful to introduce new polynomials $Q_i$ in the free cumulants. 
Define $Q_0=1$, $Q_1=0$ and for any $n \ge 2$,
\begin{equation*}
Q_n= \sum_{|\mu|=n} (l(\mu)-1)! \, \mathcal{R}_\mu.
\end{equation*}
Writing for short 
\[\mathcal{Q}_\mu=\prod_{i\ge 2} Q_i^{m_i(\mu)}/m_i(\mu)!,\quad\quad \mathcal{C}_\mu=\prod_{i\ge 2} C_i^{m_i(\mu)}/m_i(\mu)!,\]
the correspondence between these three families is given by
\begin{equation*}
\begin{split}
Q_n&= \sum_{|\mu|=n} {(-1)}^{l(\mu)} (l(\mu)-1)!  \,  \mathcal{C}_\mu,\\
C_n&= \sum_{|\mu|=n} l(\mu)! \, \mathcal{R}_\mu= \sum_{|\mu|=n}
\mathcal{Q}_\mu,\\
(1-n)R_{n}&=\sum_{|\mu|=n} {(-1)}^{l(\mu)} \,  \mathcal{Q}_\mu=
\sum_{|\mu|=n} {(-1)}^{l(\mu)} l(\mu)!\,  \mathcal{C}_\mu.
\end{split}
\end{equation*}

These relations are better understood by using the theory of symmetric functions.
Actually let $\mathbf{A}$ be the (formal) alphabet defined by 
\[(i-1)R_i=-h_i(\mathbf{A}),\quad Q_i=-p_i(\mathbf{A})/i,\quad C_i={(-1)}^i e_i(\mathbf{A}).\]
Writing 
\[u_\mu=l(\mu)!/\prod_{i\ge 1} m_i(\mu)!,\quad
 \epsilon_\mu={(-1)}^{n-l(\mu)},\quad 
 z_\mu=\prod_{i\ge 1} i^{m_i(\mu)} m_i(\mu)!,\]
the previous relations are merely the classical properties~\cite[pp. 25 and 33]{Ma}
\begin{equation*}
\begin{split}
p_n&=-n \sum_{|\mu|=n} {(-1)}^{l(\mu)} \,  u_\mu h_\mu/l(\mu)
=-n \sum_{|\mu|=n} \epsilon_\mu  u_\mu e_\mu/ l(\mu),\\
e_n&= \sum_{|\mu|=n} \epsilon_\mu u_\mu h_\mu= \sum_{|\mu|=n} \epsilon_\mu z_{\mu}^{-1} p_{\mu},\\
h_n&=\sum_{|\mu|=n} z_{\mu}^{-1} p_{\mu}=
\sum_{|\mu|=n} \epsilon_\mu u_\mu e_\mu.
\end{split}
\end{equation*}

From these relations, it is clear that $C$-positivity implies $Q$-positivity, which itself implies $R$-positivity. In particular the following conjecture is \textit{a priori} stronger than the one of Kerov-Biane and weaker than the one of Goulden-Rattan.

\begin{conj}
For $k\ge 1$ the coefficients  of $K_{r,r-2k+1}$ in terms of the $Q_i$'s are nonnegative rational numbers.
\end{conj}

The assertion is trivial for $k=1$ since
\[K_{r,r-1}= \frac{1}{4}\binom{r+1}{3}C_{r-1}=\frac{1}{4}\binom{r+1}{3}\sum_{|\mu|=r-1}
\mathcal{Q}_\mu.\]
This leads us to the following conjecture (with obviously $g_1=1/4$).
\begin{conj}
For any $k\ge 1$ there exists an inhomogeneous symmetric function $g_k$, having maximal degree $4(k-1)$, such that
\[K_{r,r-2k+1}= \binom{r+1}{3} \sum_{|\mu|=r-2k+1} (2k-1)^{l(\mu)} \, g_k(\mu)  \,  \mathcal{Q}_\mu,\]
where $g_k(\mu)$ denotes the value of $g_k$ at the integral vector $\mu$. This symmetric function is independent of $r$.
\end{conj}

It is a highly remarkable fact that, in contrast with the $C$-expansion, the $Q$-positivity is completely analogous to the $R$-positivity (and possibly equivalent).
\begin{conj}
For $k\ge 2$ the inhomogeneous symmetric function $g_k$ may be written
\[g_k= \sum_{|\rho| \le 4(k-1)} a_\rho^{(k)} m_\rho,\]
where the coefficients $a_\rho^{(k)}$ are positive rational numbers.
\end{conj}

The assertion of Conjecture 5 is a direct consequence. Conjecture 7 is supported by the following result for $k=2$, which will be proved in Section 6.

\begin{theo}
For $k=2$, we have
\begin{multline*}
8640 \,g_2=
9 m_{4}+20 m_{31}+22 m_{22}+28 m_{21^2}+24 m_{1^4}\\
+60 m_{3}+84 m_{21}+72 m_{1^3}+105 m_{2}+90 m_{1^2}+54 m_{1}.
\end{multline*}
\end{theo}

Conjecture 7 is also supported by computer calculations, giving the positive numbers $a_\rho^{(k)}$ for $k=3,4$.

\begin{conj}
For $k=3$, the values of $500 . 5! . 7 ! \, a_\rho^{(3)}$ are given by the table below.
\end{conj}
\small

\hspace{1 cm}\\
\begin{tabular}{|c|c|c|c|c|c|c|c|c|c|}
\hline
8 & 71 & 62 & $61^2$ & 53 & 521 & $51^3$ & $4^2$ & 431 & $42^2$ \\ \hline
1125 & 5400 & 13500 & 21480 & 23400 & 46200 & 69072 & 28350 &
69000 & 84900 \\ \hline
$421^2$ & $41^4$ & $3^22$ & $3^21^2$ & $32^21$ & $321^3$ & $31^5$ & $2^4$ & $2^31^2$ & $2^21^4$ \\ \hline
126168 & 174864 & 104400 & 157152 & 190704 & 265632 & 338880 & 233208 & 322128 & 414432 \\ \hline
$21^6$ & $1^8$\\ \cline{1-2}
486720 & 524160\\ \cline{1-2}
\end{tabular}

\hspace{0.4 cm}\\
\begin{tabular}{|c|c|c|c|c|c|c|c|c|}
\hline
7 & 61 & 52 & $51^2$ & 43 & 421 & $41^3$ & $3^21$ & $32^2$ \\ \hline
27000 & 107400 & 231000 & 345360 & 345000 & 630840 & 874320 & 785760 & 953520 \\ \hline
$321^2$ & $31^4$ & $2^31$ & $2^21^3$ & $21^5$ & $1^7$\\ 
\cline{1-6}
1328160 & 1694400 & 1610640 & 2072160 & 2433600 & 2620800\\
\cline{1-6}
\end{tabular}

\hspace{0.4 cm}\\
\begin{tabular}{|c|c|c|c|c|c|c|c|c|}
\hline
6 & 51 & 42 & $41^2$ & $3^2$ & 321 & $31^3$ & $2^3$ & $2^21^2$ \\ \hline
261750 & 840300 & 1532250 & 2121660 & 1907400 & 3217080 & 4095696 & 3896460 & 5001672\\ \hline
$21^4$ & $1^6$\\ \cline{1-2}
5853744 & 6274080\\ \cline{1-2} 
\end{tabular}

\hspace{0.4 cm}\\
\begin{tabular}{|c|c|c|c|c|c|c|}
\hline
5 & 41 & 32 & $31^2$ & $2^21$ & $21^3$ & $1^5$ \\ \hline
1323500 & 3322300 & 5017400 & 6358480 & 7740360 & 8988720 & 9530400 \\ \hline
\end{tabular}

\hspace{0.4 cm}\\
\begin{tabular}{|c|c|c|c|c|}
\hline
4 & 31 & $2^2$ & $21^2$ & $1^4$ \\ \hline
3763625 & 7093100 & 8590950 & 9830340 & 10212600 \\ \hline
\end{tabular}

\hspace{0.4 cm}\\
\begin{tabular}{|c|c|c|}
\hline
3 & 21 & $1^3$ \\ \hline
6011500 & 8045700 & 8043000 \\ \hline
\end{tabular}

\hspace{0.4 cm}\\
\begin{tabular}{|c|c|c|}
\hline
2 & $1^2$ & 1 \\ \hline
4985500 & 4595000 & 1650000 \\ \hline
\end{tabular}\\
\hspace{0.3 cm}\\
\normalsize

This conjecture has been checked for any $K_r$ with $r\le 32$. Starting from~\cite[Theorem 3.5.12]{R}, it may probably be proved by the method given in the next section.

We have also obtained the values of the positive numbers $a_\rho^{(4)}$. Listing them here would be tedious, but they are available upon request to the author.

\newpage
\section{Proof of Theorem 2}

We start from the following lemma of symmetric function theory. It is better understood in the language of $\lambda$-rings. This method allows to handle symmetric functions acting on ``sums'', ``products'' or ``multiples'' of alphabets. Here we shall not enter into details, and refer the reader to~\cite[Chapter 2]{Las} or~\cite[Section 3]{La1} for a short survey.

If $f$ is a symmetric function, we denote $f[\mathbf{A}]$ its $\lambda$-ring action on the alphabet $\mathbf{A}$, which should not be confused with its evaluation $f(\mathbf{A})$. For instance $p_n[-z+2]=-z^n+2$ and $p_n(-z+2)=(-z+2)^n$.

\begin{lem}
On any alphabet $\mathbf{A}$ and for any positive integer $n$, we have
\begin{equation*}
\begin{split}
\sum_{\begin{subarray}{c}(i,j,k) \in \mathsf{N}^3\\i+j+k=n
\end{subarray}} e_ie_je_k
&= \sum_{|\mu|=n} (-1)^{n-l(\mu)}\, 3^{l(\mu)}
\,z_{\mu}^{-1} p_{\mu},\\
\sum_{\begin{subarray}{c}(i,j,k) \in \mathsf{N}^3\\i+j+k=n
\end{subarray}} i^2 e_ie_je_k 
&=  \sum_{|\mu|=n} (-1)^{n-l(\mu)}\, 3^{l(\mu)-2}\, \big(n^2+2p_2(\mu)\big) \, z_{\mu}^{-1} p_{\mu}.
\end{split}
\end{equation*}
\end{lem}

\begin{proof}[Sketch of proof] Recall the ``Cauchy formula''~\cite[(1.6.6)]{Las}, or~\cite[(4.1) p. 62-65]{Ma} or~\cite[p. 222]{La1},
\begin{equation*}
e_{n} [\mathbf{AB}]= \sum_{\left|{\mu }\right| = n} 
(-1)^{n-l(\mu)}  \, z_{\mu}^{-1} p_{\mu} [\mathbf{A}]
\, p_{\mu} [\mathbf{B}].
\end{equation*}
The first relation evaluates $e_n[3\mathbf{A}]$ by using this formula together with the identity $p_{\mu}[p]=p^{l(\mu)}$ valid for any real number $p$.

For the second relation, we evaluate similarly $e_n[(z+2)\mathbf{A}]$. Then we differentiate two times and fix $z=1$. At the left-hand side we get $\sum_{i+j+k=n} i(i-1) e_ie_je_k[\mathbf{A}]$. At the right-hand side, we compute
\[\partial^2_z\big(p_\mu[z+2]\big)\Big|_{z=1}=
\partial^2_z\big(\prod_{i\ge 1} (z^i+2)^{m_i(\mu)}\big)\Big|_{z=1}
=3^{l(\mu)-2}\, \big(n^2-3n+2p_2(\mu)\big).\]
Observe that by differentiating $r$ times, we might similarly get 
$\sum_{i+j+k=n} \binom{i}{r} e_ie_je_k$.
\end{proof}

Specializing the alphabet $\mathbf{A}$ as in Section 5, so that
\[Q_i=-p_i(\mathbf{A})/i,\quad\quad C_i={(-1)}^i e_i(\mathbf{A}),\]
we obtain 
\[\sum_{\begin{subarray}{c}(i,j,k) \in \mathsf{N}^3\\i+j+k=n
\end{subarray}} (a+bi+ci^2)\, C_iC_jC_k = \sum_{|\mu|=n} 3^{l(\mu)} \, \Big(a+\frac{b}{3}n+\frac{c}{9}(n^2+2p_2(\mu))\Big)\, \mathcal{Q}_\mu.\]
By insertion in the expression~\cite[Theorem 3.3]{GR}
\[K_{r,r-3}=\binom{r+1}{3} \sum_{\begin{subarray}{c}(i,j,k) \in \mathsf{N}^3\\i+j+k=r-3\end{subarray}}
\big(a(r) +b(r) i^2\big)C_iC_jC_k,\]
we obtain
\[g_2(\mu)=a(r)+\frac{1}{9}b(r)\big((r-3)^2+2p_2(\mu)\big).\]
Since $|\mu|=p_1(\mu)=r-3$, this can be rewritten
\[g_2=\frac{1}{8640}\Big(8p_2p_1^2+p_1^4+48p_2p_1+12p_1^3+60p_2+45p_1^2+54p_1\Big).\]
Using ACE~\cite{V} the conversion to monomial functions is performed immediately.\qed

\section{Theorem 1 revisited}

In Section 3 (proof of Theorem 1) we used the property
\[\sum_{\begin{subarray}{c}(i,j,k) \in \mathsf{N}^3\\i+j+k=n
\end{subarray}} (a+bi+ci^2)\, C_iC_jC_k = \frac{1}{2} \sum_{|\mu|=n} (l(\mu)+2)! \, \Big(a+\frac{b}{3}n+\frac{c}{6}(n^2+p_2(\mu))\Big)\, \mathcal{R}_\mu,\]
which may also be proved by $\lambda$-rings method. It is obtained by specialization of the following lemma.
\begin{lem}
On any alphabet $\mathbf{A}$ and for any positive integer $n$, we have
\begin{equation*}
\begin{split}
\sum_{\begin{subarray}{c}(i,j,k) \in \mathsf{N}^3\\i+j+k=n
\end{subarray}} e_ie_je_k
&= \frac{1}{2} \sum_{|\mu|=n} (-1)^{n-l(\mu)}\, \frac{(l(\mu)+2)!}{\prod_{i} m_i(\mu)!}
\, h_\mu,\\
\sum_{\begin{subarray}{c}(i,j,k) \in \mathsf{N}^3\\i+j+k=n
\end{subarray}} i^2 e_ie_je_k 
&= \frac{1}{12} \sum_{|\mu|=n} (-1)^{n-l(\mu)}\, \frac{(l(\mu)+2)!}{\prod_{i} m_i(\mu)!} \, \big(n^2+p_2(\mu)\big) \, h_\mu.
\end{split}
\end{equation*}
\end{lem}

\begin{proof}[Sketch of proof] Recall the ``Cauchy formula''~\cite[(1.6.3)]{Las}, or~\cite[(4.2) p. 62-65]{Ma} or~\cite[p. 222]{La1},
\begin{equation*}
(-1)^n e_{n} [\mathbf{AB}]= \sum_{\left|{\mu }\right| = n} 
m_{\mu} [-\mathbf{B}] \, h_{\mu} [\mathbf{A}].
\end{equation*}
The first relation evaluates $e_n[3\mathbf{A}]$ by using this formula and the identity~\cite[(2.2.2)]{Las} valid for any real number $p$,
\[m_{\mu}[p]=p(p-1)\cdots(p-l(\mu)+1)/\prod_{i} m_i(\mu)!.\]

For the second relation, we evaluate similarly $e_n[(z+2)\mathbf{A}]$. Then we differentiate two times and fix $z=1$, obtaining $\sum_{i+j+k=n} i(i-1) e_ie_je_k[\mathbf{A}]$ at the left-hand side. At the right-hand side, we compute 
\[\prod_{i\ge 1} m_i(\mu)!\, \partial^2_z
\big(m_{\mu}[-z-2]\big)\Big|_{z=1}= (-1)^{l(\mu)} (l(\mu)+2)! \, \big(n^2-2n+p_2(\mu)\big)/12.\]
\end{proof}

\section{Final remark}
In this note, we have considered three conjectural developments of the Kerov component $K_{r,r-2k+1}$, namely up to $\binom{r+1}{3}$,
\begin{equation*}
\begin{split}
\sum_{\begin{subarray}{c}\nu \in \mathsf{N}^{2k-1}\\|\nu|=r-2k+1\end{subarray}} F_k(\nu)
\,  \prod_{i=1}^{2k-1}C_{\nu_i}&=
\sum_{|\mu|=r-2k+1} (l(\mu)+2k-2)! \, f_k(\mu) \, \mathcal{R}_\mu\\
&=\sum_{|\mu|=r-2k+1} (2k-1)^{l(\mu)} \, g_k(\mu) \, \mathcal{Q}_\mu.
\end{split}
\end{equation*}
As indicated above, these relations are better understood in the framework of symmetric functions. Choosing
\[(i-1)R_i=-h_i(\mathbf{A}),\quad Q_i=-p_i(\mathbf{A})/i,\quad C_i={(-1)}^i e_i(\mathbf{A}),\]
they are the specializations at $\mathbf{A}$ of the abstract identities 
\begin{equation*}
\begin{split}
\sum_{\begin{subarray}{c}\nu \in \mathsf{N}^{2k-1}\\|\nu|=n
\end{subarray}} F_k(\nu)
\,  e_\nu &=
\sum_{|\mu|=n} (-1)^{n-l(\mu)}\, f_k(\mu) \, \frac{(l(\mu)+2k-2)!}{\prod_{i} m_i(\mu)!} \,h_\mu \\
&=\sum_{|\mu|=n} (-1)^{n-l(\mu)}\, g_k(\mu)\,  (2k-1)^{l(\mu)} \,  z_{\mu}^{-1} p_{\mu}.
\end{split}
\end{equation*}
Moreover these identities are themselves related with the classical Cauchy formulas. Using the values of $p_\mu[p]$ and $m_\mu[p]$ given above, they may be written
\begin{equation*}
\begin{split}
 (-1)^n \sum_{\begin{subarray}{c}\nu \in \mathsf{N}^{2k-1}\\
|\nu|=n\end{subarray}} F_k(\nu)
\,  e_\nu &=\sum_{|\mu|=n} f_k(\mu) \, m_\mu[-2k+1] \, h_\mu\\
&=\sum_{|\mu|=n} g_k(\mu) \, p_\mu[-2k+1] \,  z_{\mu}^{-1} p_{\mu}.
\end{split}
\end{equation*}

Therefore it seems plausible that the conjectured positivity properties of $f_k$ and $g_k$ are equivalent, and reflect some abstract pattern of the theory of symmetric functions.


\begin{thebibliography}{29}
\bibitem{B1}
P.\ Biane, \emph{Representations of symmetric groups and free probability}, Adv. Math. \textbf{138} (1998), 126Ð-181.
\bibitem{B2}
P.\ Biane, \emph{Characters of symmetric groups and free cumulants}, Lecture Notes in Math. \textbf{1815} (2003), 185Ð-200, Springer, Berlin, 2003. 
\bibitem{B3}
P.\ Biane, \emph{On the formula of Goulden and Rattan for Kerov polynomials}, S\'em. Lothar. Combin., \textbf{55} (2006), article B55d.
\bibitem{F}
V.\ F\'eray, \emph{Combinatorial interpretation and positivity of Kerov's character polynomials}, arXiv 0710.5885.
\bibitem{GR}
I.\ P.\ Goulden, A.\ Rattan, \emph{An explicit form for Kerov's character polynomials}, Trans. Amer. Math. Soc. \textbf{359} (2007), 3669--3685. 
\bibitem{K}
S.\ V.\ Kerov, talk at IHP Conference (2000).
\bibitem{Las}
A.\ Lascoux, \emph{Symmetric functions and combinatorial operators on
polynomials}, CBMS Regional Conference Series in Mathematics \textbf{99}, Amer. Math. Soc., Providence, 2003.
\bibitem{La1}
M.\ Lassalle, \emph{Une $q$-sp\'ecialisation pour les fonctions
sym\'etriques monomiales}, Adv.\ Math. \textbf{162} (2001), 217--242.
\bibitem{Ma}
I.\ G.\ Macdonald, \emph{Symmetric Functions and Hall Polynomials},
Clarendon Press, second edition, Oxford, 1995.
\bibitem{R}
A.\ Rattan, \emph{Character polynomials and Lagrange inversion}, Thesis (2005), Waterloo University. 
\bibitem{S}
P.\ \.{S}niady, \emph{Asymptotics of characters of symmetric groups and free probability}, Discrete Math. \textbf{306} (2006), 624--665.
\bibitem{V}
S.\ Veigneau, \emph{ACE, an Algebraic Combinatorics Environment for the computer algebra system Maple}, available at http://phalanstere.univ-mlv.fr/{\textasciitilde}ace/

\end{thebibliography}
\end{document}